\documentclass[11pt]{article}
\newcommand{\be}{\begin{equation}}
\newcommand{\ee}{\end{equation}}
\newcommand{\bea}{\begin{eqnarray}}
\newcommand{\eea}{\end{eqnarray}}
\newcommand{\bean}{\begin{eqnarray*}}
\newcommand{\eean}{\end{eqnarray*}}
\newcommand{\brray}{\begin{array}}
\newcommand{\erray}{\end{array}}
\newcommand{\ben}{\begin{equation}{nonumber}}
\newcommand{\een}{\end{equation}{nonumber}}
\newcommand{\newsection}[1]{\setcounter{dfn}{0}
\section{#1}}

\newtheorem{dfn}{Definition}[subsection]
\newtheorem{thm}[dfn]{Theorem}
\newtheorem{lmma}[dfn]{Lemma}
\newtheorem{ppsn}[dfn]{Proposition}
\newtheorem{crlre}[dfn]{Corollary}
\newtheorem{xmpl}[dfn]{Example}
\newtheorem{rmrk}[dfn]{Remark}

\newcommand{\bdfn}{\begin{dfn}}
\newcommand{\bthm}{\begin{thm}}
\newcommand{\blmma}{\begin{lmma}}
\newcommand{\bppsn}{\begin{ppsn}}
\newcommand{\bcrlre}{\begin{crlre}}
\newcommand{\bxmpl}{\begin{xmpl}}
\newcommand{\brmrk}{\begin{rmrk}}

\newcommand{\edfn}{\end{dfn}}
\newcommand{\ethm}{\end{thm}}
\newcommand{\elmma}{\end{lmma}}
\newcommand{\eppsn}{\end{ppsn}}
\newcommand{\ecrlre}{\end{crlre}}
\newcommand{\exmpl}{\end{xmpl}}
\newcommand{\ermrk}{\end{rmrk}}

\newcommand{\IC}{{C\!\!\!|~}}

\newcommand{\IZ}{{Z \! \!\!\!Z}}

\newcommand{\al}{\alpha}

\newcommand{\cla}{{\cal A}}
\newcommand{\clb}{{\cal B}}
\newcommand{\clc}{{\cal C}}

\newcommand{\cle}{{\cal E}}
\newcommand{\clf}{{\cal F}}

\newcommand{\clh}{{\cal H}}

\newcommand{\clk}{{\cal K}}
\newcommand{\cll}{{\cal L}}
\newcommand{\clm}{{\cal M}}

\newcommand{\clt}{{\cal T}}

\def\a*{{\cal A}_{h,*}}
\def\B{{\cal B}(h)}
\def\B1{{\cal B}_1(h)}
\def\b{{\cal B}^{s. a. }(h)}
\def\b1{{\cal B}^{s. a. }_1(h)}

\def\a0{{\cal A}_0}
\def\c0{{\cal C}_0}

\newcommand{\ot}{\otimes}

\newcommand{\raro}{\rightarrow}

\newcommand{\id}{\mbox{id}}

\newcommand{\qed} { \mbox{}\hfill \vspace{1ex}}

\begin{document}
 \begin{center}
{\large {\bf A Complete Formulation of Baum-Connes' Conjecture for
the Action of Discrete Quantum Groups }}\\ by\\ {\large Debashish
Goswami}\\ {Stat-Math Unit, Indian Statistical Institute,}\\ {203,
B.T. Road, Kolkata 700 108, India}\\{\large e-mail :
goswamid@isical.ac.in}\\ and\\ {\large A. O. Kuku;}\\ {\large The
Abdus Salam International Centre for Theoretical Physics
(Mathematics Section)}\\ {\large Strada Costiera 11, trieste
34014, Italy.}\\  {\large
 e-mail : kuku@ictp.trieste.it}\\
\end{center}
\begin{abstract}
 We formulate a version of Baum-Connes' conjecture for a discrete
 quantum group, building on our earlier work (\cite{GK}).  Given such a quantum group $\cla$, we construct
 a directed family $\{ \cle_F \}$ of $C^*$-algebras ($F$ varying over some suitable
 index set), borrowing the ideas of \cite{cuntz}, such that there
 is a natural action of $\cla$ on each $\cle_F$ satisfying the
 assumptions of \cite{GK}, which makes it possible to define the
 ``analytical assembly map", say $\mu^{r,F}_i$, $i=0,1,$  as in \cite{GK},
  from the $\cla$-equivariant $K$-homolgy groups of $\cle_F$ to the $K$-theory groups
   of the ``reduced" dual $\hat{\cla_r}$ (c.f. \cite{GK} and the references therein for
    more details).   As a result, we can define the Baum-Connes' maps
    $\mu^r_i
   : \stackrel{\rm lim}{\longrightarrow} KK_i^\cla(\cle_F,\IC) \raro
   K_i(\hat{\cla_r}),$ and in the classical case, i.e. when
   $\cla$ is $C_0(G)$ for a discrete group, the isomorphism of the
   above maps for $i=0,1$ is equivalent to the Baum-Connes'
   conjecture. Furthermore, we verify its truth for an arbitrary
   finite dimensional quantum group and obtain partial results for
   the dual of $SU_q(2).$

\end{abstract}
\vspace{7mm} {\large Key words : Baum-Connes Conjecture, Discrete
Quantum Group,
Equivariant KK-Theory}

  \section{Introduction}
\subsection{Classical Baum-Connes' Conjecture}
Baum-Connes' conjecture has been occupying the centre-stage of
K-Theory and geometry for more than two decades.  Given a locally
compact group $G$, and a locally compact Hausdorff space $X$
equipped with a $G$-action such that $X$ is proper and
$G$-compact, there are canonical maps $\mu^r_i :
KK_i^G(C_0(X),\IC) \raro K_i(C^*_r(G))$, and $\mu_i :
KK_i^G(C_0(X),\IC) \raro K_i(C^*(G)),$ for $i=0,1,$ where $C_0(X)$
is the commutative $C^*$-algebra of continuous complex-valued
functions on $X$ vanishing at infinity, $C^*_r(G)$ and $C^*(G)$
are respectively the reduced and free groups $C^*$-algebras, and
$KK^G_.$ denotes the Kasparov's equivariant KK-functor. In
particular, $KK^G(C_0(X),\IC)$ is identified with the
$G$-equivariant K-homology of $X$, and thus is essentially
something geometric or topological, whereas the object
 $K_i(.)$ on the right hand side involves the reduced or free group algebras,
which are analytic in some sense.

Now, let $\underline{EG}$ be the universal space for proper
actions of $G$. The definition of proper $G$-actions  and explicit
constructions in various cases of interest can be found in
\cite{val1},\cite{val2} and the references therein. The
equivariant $K$-homology of $\underline{EG}$, say
$RK^G_i(\underline{EG}),i=0,1,$ can be defined as the inductive
limit of $KK^G_i(C_0(X),\IC),$ over all possible locally compact,
$G$-proper and $G$-compact subsets $X$ of the universal space
$\underline{EG}.$  Since the construction of $KK^G_i$ and $K_i$
commute with the procedure of taking an inductive limit, it is
possible to define $\mu^r_i,\mu^i$ on the equivariant $K$-homology
$RK^G_i(\underline{EG})$.   The conjecture of Baum-Connes states
that $\mu^r_i,i=0,1$ are isomorphisms of abelian groups. This
conjecture admits certain other generalizations, such as the
Baum-Connes conjecture with coefficients (which seems to be false
from some recent result announced by M. Gromov, see \cite{val1}
for references), but we do not want to discuss those here.
However, we would like to point out that the Baum-Connes
conjecture has already been verified for many classical groups
(and for a wide variety of coefficient algebras as well) , using
different methods and ideas from many diverse areas of
mathematics, and has given birth to many new and interesting tools
and techniques in all these areas. In fact, the truth of this
conjecture, if established, will prove many other famous
conjectures in topology, geometry and $K$-theory.

\subsection{Motivation for a quantum version}
 During the  last two decades, the theory of quantum groups, which is a
natural and far-reaching generalization of the concept of
topological groups, has become another fast-growing branch of
mathematics and mathematical physics, thanks to the works by
Woronowicz, Drinfeld, Jimbo and many other mathematicians
(\cite{wor}, \cite{drin}, \cite{jim}). On the other hand, with the
pioneering efforts of Connes (see \cite{Con}), followed by himself
and many other mathematicians, a powerful generalization of
classical differential and Riemannian geometry has emerged under
the name of noncommutative geometry, which has had, since its very
beginning, very close connections with $K$-theory too.
Furthermore, Baaj and Skandalis (\cite{bs}) have been able to
construct an analogue of equivariant $KK$-theory for the actions
of quantum groups,
 as natural extension of Kasparov's equivariant $KK$-theory. This motivates
one to think of a possibility of generalizing the Baum-Connes
construction in the framework of quantum groups, as one hopes for
a possible rich interplay between noncommutative geometry and
quantum groups in this context. In the present article, we make an
attempt towards this generalization, completing our formulation of
this conjecture for discrete quantum groups as a follow-up of our
earlier paper \cite{GK} in this direction.

\subsection{ Plan of the paper}

 As remarked in \cite{GK}, the classical formulation of the
 Baum-Connes conjecture for the action of locally compact groups
 or discrete groups $G$ could be achieved in two steps. The first
 step is to define assembly maps $\mu_i, \mu_i^r$ for $G$-compact
 and $G$-proper actions. We achieved the analogue of this step for
 discrete quantum groups in \cite{GK}. More precisely, given an
 action of a discrete quantum group $\cla$ on a $C^*$-algebra
 $\clc$ satisfying certain regularity assumptions (resembling the
 notion of proper $G$-compact action of classical discrete groups
 on some space), we constructed canonical maps $\mu_i,\mu_i^r$
 ($i=0,1$) from the $\cla$-equivariant $K$-homology groups
 $KK_i^\cla(\clc, \IC)$ to the $K$-theory groups
 $K_i(\hat{\cla}),K_i(\hat{\cla}_r)$ respectively, where
 $\hat{\cla},\hat{\cla_r}$ denote respectively the quantum
 analogue of the full and reduced group $C^*$-algebras.
 For readers' convenience, we briefly review in section 4 our
 constructions and results in \cite{GK}. We also illustrate these
  constructions with examples in 4.3.

 The second step in the classical formulation of the Baum-Connes
 conjecture is to define a universal space $\underline{EG}$ for
 the proper actions of $G$ and build explicit good models for
 $\underline{EG}$ to show that it can be approximated by its
 subsets $X$ having $G$-proper and $G$-compact actions, thereby
 defining appropriate maps from the equivariant $K$-homology
 groups $RK_i^G(\underline{EG})$ = $\stackrel{\rm lim}{\longrightarrow}
 KK_i^G(C_0(X),\IC)$ to
  $K_i(\hat{\cla}), K_i(\hat{\cla}_r)$ by inductive limits. We
  achieve this step in this paper for the action of discrete
  quantum groups by borrowing ideas from the work of Cuntz
  \cite{cuntz}, who used the language of noncommutative simplicial
  complexes to describe this step in the classical situation. More
  precisely, we construct a directed family $\{ \cle_F \}$ of
  $C^*$-algebras ($F$ varying over some index set) and prove that
  the natural action of $\cla$ on $\cle_F$ satisfies the
  assumptions in \cite{GK} which makes it possible to define
  analytic assembly maps $\mu_i^F, \mu_i^{r,F} (i=0,1)$ as in
  \cite{GK} from the $\cla$-equivariant $K$-homology groups of
  $\cle_F$ to the $K$-theory groups of $\hat{\cla}$ and
  $\hat{\cla}_r$. Consequently, we are able to define the
  Baum-Connes maps
  $$ \mu_i : \stackrel{\rm lim}{\longrightarrow} KK_i^\cla(\cle_F,\IC) \raro
  KK_i(\IC,\hat{\cla}),$$
  $$ \mu_i^r : \stackrel{\rm lim}{\longrightarrow}
  KK_i^\cla(\cle_F,\IC) \raro KK_i(\IC,\hat{\cla}_r);$$
   so that in the classical case when $\cla=C_0(G),$ the
   isomorphism of $\mu_i^r$ is equivalent to the Baum-Connes
   conjecture (see Proposition 5.2.1).

   We verify the conjecture for an arbitrary finite dimensional
   quantum groups (Theorem 5.2.2),
     and also give some interesting examples and computations of the analytical assembly map (see 4.3 and Remark 5.2.3).

\section{Notes on notation}
We shall mostly follow the notation used in \cite{GK}. However,
for readers' convenience, let us briefly mention some of those
here.

 \subsection{Hilbert spaces and operators :} All the Hilbert spaces considered
  in this paper are over the field of complex numbers.

For two Hilbert spaces $\clh_1,\clh_2$ and some bounded operator
$X \in \clb(\clh_1 \ot \clh_2)=\clb(\clh_1) \ot \clb(\clh_2),$ we
denote by $X_{12}$ the operator $ X \ot 1_{\clh_2}$ on $\clh_1 \ot
\clh_2 \ot \clh_2,$ and denote by $X_{13}$ the operator
$(1_{\clh_1} \ot \Sigma) (X \ot 1_{\clh_2})(1_{\clh_1} \ot
\Sigma)$ on $\clh_1 \ot \clh_2 \ot \clh_2$, where $\Sigma : \clh_2
\ot \clh_2 \raro \clh_2 \ot \clh_2$ flips the two copies of
$\clh_2.$ For two vectors $\xi,\eta \in \clh_1$ we define a map
$T_{\xi \eta} : \clb(\clh_1 \ot \clh_2) \raro \clb(\clh_2)$ by
setting $T_{\xi \eta}(A \ot B):=<\xi,A\eta>B,$ where $A \in
\clb(\clh_1), B \in \clb(\clh_2),$ and extend this definition to
the whole of $\clb(\clh_1 \ot \clh_2)$ in the obvious way.
 For some  Hilbert space $\clh$, we
denote by  $\clb_0(\clh)$  the $C^*$-algebra of compact operators
on $\clh$.

\subsection{Multiplier algebras}
 For a Hilbert space $\clh$, and a pre-$C^*$-algebra
$\a0 \subseteq \clb(\clh),$ we shall denote the multiplier algebra
of the norm-closure of $\a0$ by $\clm(\a0).$ $\clm(\a0)$ has two
natural topologies : one is the norm topology in which it becomes
a unital (but typically nonseparable) $C^*$-algebra, and the other
is the strict topology, which makes it  a locally convex
 topological $\ast$-algebra.

 The case when $\a0$ is an algebraic direct sum of finite
 dimensional matrix algebras is of special interest to us. For
 such an algebra, say of the form, $\a0=\oplus_{\al \in I}
 \cla_\al$, where $I$ in some index set and $\cla_\al$ is $n_\al
 \times n_\al$ complex matrix algebra ($n_\al $ positive integer),
  the multiplier algebra $\clm(\a0)$ can be described as the set
  of
 all collections
$(a_\al)_{\al \in I}$ with $a_\al \in \cla_\al$ for each $\al$,
and $\sup_\al \| a_\al \| < \infty.$ The $\ast$-algebra operations
are taken to be the obvious ones; i.e.
$(a_\al)+(b_\al):=(a_\al+b_\al)$, $(a_\al).(b_\al):=(a_\al b_\al)$
and $(a_\al)^*:=(a_\al^*).$   Similarly,  $\clm(\a0 \ot \a0)$
consists  of all collections of the form $(a_\al \ot b_\beta)$
where $\al, \beta $ vary  over $I.$

\subsection{Algebraic multiplier}
We shall also use the following algebraic analogue of multiplier
 algebra.
 \bdfn
 An ``algebraic multiplier" of the $\ast$-algebra $\a0$ above is
  a collection $(a_\al)_{\al \in I},$ with $a_\al
\in \cla_\al \forall \al$ (no restriction on norms).
 The set of all algebraic multipliers of $\a0$ will be  denoted  by $\clm_{\rm alg}(\a0)$, and
 we equip it with the structure of a
 $\ast$-algebra by defining the following obvious operations :
$$
{\rm For} \hspace{2mm} U =(u_\al),V=(v_\al) \in \clm_{\rm
alg}(\a0),\hspace{2mm} UV:=(u_\al v_\al)_\al,\hspace{2mm}
U^*:=(u_\al^*).$$ \edfn

\brmrk
  Clearly, any
element of $\a0$ can be viewed as an element of $\clm_{\rm
alg}(\a0)$, by thinking of $a \in \a0$ as $(a_\al)$, where $a_\al$
is the component of $a$ in $\cla_\al.$ It is easy to see that
$Ua,aU \in \a0$ for $a \in \a0, U \in \clm_{\rm alg}(\a0).$
Similarly,  an algebraic multiplier of $\a0 \ot \a0$  is given by
any collection of the form $M \equiv (m_{\al \beta})_{\al,\beta
\in I}$, with $m_{\al \beta} \in \cla_\al \ot \cla_\beta.$
 \ermrk
\brmrk Let $\clk$ be the smallest Hilbert space containing the
algebraic direct sum $\oplus_{\al \in I} \clk_\al \equiv
\oplus_\al \IC^{n_\al},$ i.e. $\clk=\{ (f_\al)_{\al \in I} : f_\al
\in \clk_\al =\IC^{n_\al}, \sum_\al \| f_\al \|^2 < \infty \},$
where the possibly uncountable sum $\sum_\al$ means the limit over
the net consisting of all possible sums over finite subsets of
$I$. Let us consider the canonical imbedding of $\a0$ in
$\clb(\clk)$, with $\cla_\al$ acting on $\IC^{n_\al}$.  Let us fix
some  matrix units $e^\al_{ij}, i,j =1,...,n_\al$ for
$\cla_\al=M_{n_\al}$, w.r.t. some fixed orthonormal basis
$e^\al_i, i=1,...,n_\al,$ of $\IC^{n_\al},$ and thus the
norm-closure of $\a0$, to be denoted by $\cla$,   is the
$C^*$-algebra generated by $e^\al_{ij}$'s. It is  easy to see that
any element of $\clm_{\rm alg}(\a0)$ can be viewed as a possibly
unbounded operator on $\clk$, with the domain containing the
algebraic direct sum of $\clk_\al$'s. Similarly, elements of
$\clm_{\rm alg}(\a0 \ot \a0)$ can be thought of as possibly
unbounded operators on $\clk \ot \clk$ with suitable domain.
\ermrk

\subsection{$C^*$ and von Neumann Hilbert modules}

 For a Hilbert $C^*$-module $E$, we denote by  $\cll(E)$  the $C^*$-algebra of adjointable
linear maps on  $E.$
 Furthermore, for a von
Neumann algebra $\clb \subseteq \clb(\clh),$ and some Hilbert
space $\clh^\prime,$  We denote by $\clh^\prime \ot \clb$ the
Hilbert von Neumann module obtained from the algebraic
$\clb$-module  $\clh^\prime \ot_{\rm alg} \clb $ by completing
this algebraic module in the strong operator topology inherited
from $\clb(\clh,\clh^\prime \ot \clh),$  where we have identified
an element of the form $(\xi \ot b)$, $\xi \in \clh^\prime, b \in
\clb,$ with the operator which sends a vector $v \in \clh$ to
$(\xi \ot bv) \in \clh^\prime \ot \clh.$

We also  introduce the following notation :
$$
{\rm For} \hspace{2mm}  \eta \in \clh^\prime, X \in
\clb(\clh^\prime) \ot \clb \equiv \cll(\clh^\prime \ot \clb),
\hspace{2mm} X \eta := X(\eta \ot 1_\clb) \in \clh^\prime \ot
\clb.$$

 Similarly, for a possibly nonunital $C^*$-algebra $\cla$, we can complete the
algebraic pre-Hilbert $\cla$ module $\clh^\prime \ot_{\rm alg}
\cla$ in the locally convex topology coming from the strict
topology on $\clm(\cla)$, so that the completion becomes in a
natural way a locally convex Hilbert $\clm(\cla)$-module, to be
denoted by $\clh^\prime \ot \clm(\cla).$ It is also easy to see
that if $X \in \clm(\clb_0(\clh^\prime) \ot \cla) ,$ $\eta \in
\clh^\prime,$ then we have $X\eta \equiv X(\eta \ot 1) \in
\clh^\prime \ot \clm(\cla).$

\subsection{Strictly continuous extensions}
If $\clb_1,\clb_2$ are two von Neumann algebras, $\clh^\prime$  a
Hilbert space,  and $\rho : \clb_1 \raro \clb_2$  a normal
$\ast$-homomorphism, then it is easy to show that $(id \ot \rho )
: \clh^\prime \ot_{\rm alg} \clb_1 \raro \clh^\prime \ot_{\rm alg}
\clb_2$ admits a unique extension (to be denoted again by $(id \ot
\rho)$) from the Hilbert von Neumann module $\clh^\prime \ot
\clb_1$ to the Hilbert von Neumann module $\clh^\prime \ot
\clb_2.$ Furthermore, one has that $(id \ot \rho)(X \eta)=(id \ot
\rho)(X) \eta$ for $X \in \clb(\clh^\prime) \ot \clb,$ $\eta \in
\clh^\prime.$ By very similar arguments one can also prove that if
$\cla_1,\cla_2$ are two $C^*$-algebras, and $\pi : \cla_1 \raro
\cla_2$ is a nondegenerate $\ast$-homomorphism (hence extends
uniquely as a unital strictly continuous $\ast$-homomorphism from
$\clm(\cla_1)$ to $\clm(\cla_2)$), then $(id \ot \pi) :
\clh^\prime \ot_{\rm alg} \cla_1 \raro \clh^\prime \ot_{\rm alg}
\cla_2$ admits a unique  extension (to be denoted by the same
notation) from $\clh^\prime \ot \clm(\cla_1)$ to $\clh^\prime \ot
\clm(\cla_2),$ which is continuous in the locally convex
topologies coming from the respective strict topologies. We also
have that $(id \ot \pi)(X \eta)=(id \ot \pi)(X) \eta,$ for $X \in
\clm(\clb_0(\clh^\prime) \ot \cla), \eta \in \clh^\prime.$

\section{Preliminaries on discrete quantum groups}

\subsection{Definition}
We briefly discuss the theory of {\it discrete quantum groups} as
developed in \cite{van}, \cite{er},\cite{van3}, \cite{kus} and
other relevant references to be found there. As in the previous
section, let us fix an index set $I$ (possibly uncountable), and
let $\a0:=\oplus_{\al \in I} \cla_\al$ be the {\it algebraic}
direct sum of $\cla_\al$'s, where for each $\al,$
$\cla_\al=M_{n_\al}$ is the finite dimensional $C^*$-algebra of
$n_\al \times n_\al$ matrices with complex entries, and $n_\al$ is
some positive integer. We also take $\clk=\oplus_\al \IC^{n_\al}$
as in the previous section.

\bdfn{{\bf Discrete Quantum Group}}

We say that $\cla$ (which is the norm-closure of $\a0$) is a
discrete quantum group if  there is a unital $C^*$-homomorphism
$\Delta : \clm(\a0) \raro \clm(\a0 \ot \a0)$ which satisfies the
following :\\  (i) For $a , b \in \a0,$ we have $$ T_1(a \ot b)
:=\Delta(a)(1 \ot b) \in \a0 \ot_{\rm alg} \a0,$$  and
$$ T_2(a \ot b) :=(a \ot 1) \Delta(b) \in \a0 \ot_{\rm alg} \a0;$$
(ii) $T_1, T_2 : \a0 \ot_{\rm alg} \a0 \raro \a0 \ot_{\rm alg}
\a0$ are
bijections;\\
(iii) $\Delta$ is coassociative in the sense that
$$ (a \ot 1 \ot 1)(\Delta \ot id)(\Delta(b)(1 \ot c))=(id \ot \Delta)((a \ot
1)\Delta(b))(1 \ot 1 \ot c),$$
 for $a,b, c \in \a0.$

 \edfn

\brmrk
 As explained in the relevant references mentioned above,
$(\Delta \ot id),$
 $(id \ot \Delta)$
admit  extensions as  $C^*$-homomorphisms from $\clm(\a0 \ot \a0
)$ to $\clm(\a0 \ot \a0 \ot \a0)$ (we denote these extensions by
the same notation) and the condition (iii) translates into
$(\Delta \ot id) \Delta=(id \ot \Delta) \Delta.$

\ermrk

 Let us recall (without proof) from \cite{van} and \cite{kus} some
of the important properties of our discrete quantum group $\cla.$
It is remarkable that it is possible to deduce from (i) to (iii)
the existence of a canonical antipode $S :\a0 \raro \a0$
satisfying $S(S(a)^*)^*=a$ and other usual properties of the
antipode of a Hopf algebra. Furthermore, there exists a counit
$\epsilon :\a0 \raro \IC.$ For details of the constructions of
these maps and their properties we refer to \cite{van}.

\blmma{ (\cite{van})}

\label{bijection}

 (i) There is a bijection of the index set
$I$, say $\al \mapsto \al^\prime$, such that
$S(e_\al)=e_{\al^\prime}, S(e_{\al^\prime})=e_\al.$

(ii) For fixed $\al,\beta \in I$, there is a finite number of
$\gamma \in I$ such that $\Delta(e_\gamma)(e_\al \ot e_\beta)$ is
nonzero.

\elmma

This allows us to make the following

 \bdfn
  Define $S(X)$ for $X =(x_\al)_I \in \clm_{\rm alg}(\a0)$ by
$S(X):=X^\prime=(x^\prime_\al)_I$ where
$x^\prime_\al=S(x_{\al^\prime}).$ Similarly, using (ii) in the
Lemma \ref{bijection},
  $\Delta(X)$ is defined as the element $Y \in
\clm_{\rm alg}(\a0 \ot \a0)$ such that $Y=(y_{\al \beta}),$ where
$y_{\al \beta}= \sum_\gamma \Delta(x_\gamma)
\Delta(e_\gamma)(e_\al \ot e_\beta).$

\edfn

\brmrk

 For algebraic multipliers
$A,B$ of $ \a0$ and $L$ of $\a0 \ot \a0,$ it is clear that
$\Delta(A)=L$ if and only if $\Delta(Aa)=L\Delta(a)$ $\forall a
\in \a0,$ and  $S(A)=B$ if and only if $S(Aa)=S(a)B$ for $a \in
\a0.$

\ermrk

\subsection{Invariant functionals and modular operator}
 Let us denote by $\a0^\prime$ the set of all linear functionals
on $\a0$ having ``finite support", i.e. they vanish on
$\cla_\al$'s for all but finite many $\al \in I.$ It is clear that
any $f \in \a0^\prime$ can be identified as a functional on
$\clm_{\rm alg}(\a0)$, by defining $f((a_\al)_I):=\sum_{\al \in I}
f(a_\al) \equiv \sum_{I_0} f(a_\al),$ where $I_0$ is the finite
set of $\al$'s such that for $\al$'s not belonging to $I_0$,
$f|_{\cla_\al}=0.$ With this identification, $f(1)$ makes sense
for any functional $f$  on $ \clm_{\rm alg}(\a0)$. Let us denote
by $e_\al$ the identity of $\cla_\al=M_{n_\al}$, which is a
minimal central projection in $\a0.$ For any subset $I_1$ of $I$
we denote by $e_{I_1}$ the direct sum of $e_\al$'s for $\al \in
I_1.$ It is clear that a functional $f$ on $\a0$ is in
$\a0^\prime$ if and only if there is some finite $I_1$ such that
$f(a)=f(e_{I_1}a)$ for all $a \in \a0.$

\bdfn

 We say that a linear functional $\phi$ (not necessarily with
finite support) on $\a0$ is left invariant if we have $(id \ot
\phi)((b \ot 1)\Delta(a))=b\phi(a)$ for all $a,b \in \a0,$ or
equivalently, $\phi((\omega \ot id)(\Delta(a)))=\omega(1) \phi(a)$
for all $a \in \a0,$ $\omega \in \a0^\prime.$ Similarly, a linear
functional $\psi$ on $\a0$ is called right invariant if $(\psi \ot
id)((1 \ot b)\Delta(a))=\psi(a)b$ for all $a,b \in \a0.$

\edfn

We now recall some of the main results regarding left and right
invariant functionals as proved in \cite{van}.

\bppsn
 Up to constant multiples,
there is a unique left invariant functional, as well as a unique
right invariant functional. However, in general (unless $S^2=id$)
left and right invariant functionals are not the same.

\eppsn

 Let us summarize some of the important and useful facts here. For
 details, we refer to \cite{van}, \cite{kus} and \cite{GK}.

 \bppsn

  (a) There exists a positive  invertible
 element $\theta \in \clm_{\rm alg}(\a0)$, identified as a possibly unbounded operator  on $\clk$,
with its domain containing all $\clk_\al$'s, such that
$\theta_\al=\theta|_{\clk_\al} \in \cla_\al$ satisfies
$\Delta(\theta)=(\theta \ot \theta),$ $S(\theta)=\theta^{-1},$ and
$S(\theta^{-1})=\theta.$\\ (b) $S^2(a)=\theta^{-1}a\theta$ for all
$a \in \a0.$\\ (c) We can choose a positive faithful left
invariant functional (to be referred to as left haar measure later
on) $\phi$ and a positive faithful right invariant functional (to
be referred to as right haar measure) $\psi$ such that
$\psi(a)=\phi(a \theta^2)=\phi(\theta^2 a)$ for $a \in \a0.$\\ (d)
$\phi(S^2(a))=\phi(a), \psi(S^2(a))=\psi(a)$ for all $a \in \a0$,
where $\phi,\psi$ as in (c).

\eppsn

 Let us now extend the definition of $\phi$ and $\psi$ on a
larger set than $\a0$ as follows.

\bdfn
 For a nonnegative element $a \in \clm(\cla_0) \subseteq
\clb(\clk),$ we define $\phi(a)$ as the limit of $\phi_J(a)$,
whenever this limit exists as a finite number, and where $J$ is
any finite subset of $I$, $\phi_J(.):=\phi(e_J.)=\phi(.e_J),$ and
the limit is taken over the net of finite subsets of $I$ partially
ordered by inclusion. Similarly, we set $\psi(a)=\lim_J
\psi(e_Ja)$ whenever the limit exists as a finite number. Since  a
general element $a \in \clm(\cla_0)$ can be    canonically written
as a linear combination of four nonnegative elements, we can
extend the definition of $\phi$ on $\clm(\cla_0)$ by linearity.
For any nonnegative $X \in \clm(\clb_0(\clh) \ot \cla)$ (where
$\clh$ is some Hilbert space), we define $(id \ot \phi)(X)$ as the
limit in the weak-operator topology  (if it exists as a bounded
operator) of the net $(id \ot \phi_J)(X)$ over finite subsets $J
\subseteq I,$ and extend this definition for a general $X \in
\clm(\clb_0(\clh) \ot \cla)$ in the usual way.    Similar
definition is given for $(id \ot \psi)$.

\edfn

 \blmma
 \cite{GK}
 If we choose $\clh=\clk$ in the
above, and take any $a \in \clm(\cla_0)$ such that $\phi(a)$ is
finite, then $(id \ot \phi)(\Delta(a))=\phi(a)1_{\clm(\cla)}.$
\elmma

We remark that an analogous fact is true for $\psi.$

The following fact, proved in  \cite{van}, will be useful later
on.

\bppsn

\label{h}

  There is a one-dimensional
 component $M_{\al_0}$ for some $\al_0 \in I$ such that the
 identity $h_0$, say, of this component has the property that
 $h_0.a=a.h_0=\epsilon(a)h_0 \forall a \in \cla_0,$ and also $\phi(h_0)=1$.
\eppsn

\subsection{Representation and dual of a discrete quantum group}
\bdfn We say that a unitary element in $\clm(\cll(\clh \ot \cla))
\equiv \clm(\clb_0(\clh) \ot \cla)$ is a unitary representation of
the discrete quantum group $\cla$ if $(id \ot \Delta)(U)=U_{12}
U_{13},$ and $(id \ot S)(U)=U^*.$ \edfn

\brmrk

 Note that the second equality in the above definition
has to be understood in the sense
 of the definition of $S$ on the algebraic multiplier, i.e. $(id \ot
S)(U(1 \ot a))=(1 \ot S(a))U^*$ for all $a \in \a0.$ Let us also
make the following useful observation : for $X \in
\clm(\clb_0(\clh) \ot \cla),$ and $\xi, \eta \in \clh,$ we have
that $T_{\xi,\eta}(X) \in \clm(\cla).$

\ermrk

 We shall now define a $\ast$-algebra structure on $\a0^\prime$,
and  identify $\a0$ with suitable elements of $\a0^\prime,$
thereby equipping $\a0$ with this new $\ast$-algebra structure,
and finally consider suitable $C^*$-completions.  This will give
rise to the analogues of the full and reduced group $C^*$-algebra
in the framework of discrete quantum groups.

\bdfn

  Define $f \ast g
$ for $f,g \in \a0$ by $(f \ast g)(a):=(f \ot g)(\Delta(a)), a \in
\a0.$ We also define an adjoint by $f^*(a):=\bar{f}(S(a)^*), a \in
\a0.$

\edfn

\brmrk

 Note that since $f,g$ have finite supports, there is some
finite subset $J$ of $I$ such that $(f \ot g)(\Delta(a))=(f \ot
g)((e_J \ot e_J)\Delta(a)),$ and since $(e_J \ot e_J)\Delta(a) \in
\a0 \ot_{\rm alg} \a0,$ $f \ast g$ is well defined.

\ermrk

 We  define for each $a \in \a0$, an element $\psi_a \in
\a0^\prime$ by $\psi_a(b):=\psi(ab).$
 It is easy to verify the following by using standard formulae involving
$\Delta$ and $S.$

\bppsn {\cite{GK}}

 For $a,b \in \a0,$ $\psi_a \ast \psi_b=\psi_{a \ast b},$ where
$a \ast b :=(id \ot \psi)((1 \ot b)((id \ot
S^{-1})(\Delta(a))))=(\phi \ot id)((a \ot 1)((S \ot
id)(\Delta(b)))).$ Furthermore, $\psi_a^*=\psi_{a^\sharp},$ where
$a^\sharp:=\theta^{-2}S^{-1}(a^*).$

\eppsn

We denote by $\hat{\a0}$ the set $\a0$ equipped with the
$\ast$-algebra structure given by $(a,b) \mapsto a \ast b, a
\mapsto a^\sharp$ described by the above proposition. There are
two different natural ways of making $\hat{\a0}$ into a
$C^*$-algebra, and thus we obtain the so-called reduced
$C^*$-algebra $\hat{\cla}_r$ and the free or full $C^*$-algebra
$\hat{\cla}$. This is done in a similar way as in the classical
case : one can realize elements of $\hat{\a0}$ as bounded linear
operators on the Hilbert space $L^2(\phi)$ (the GNS-space
associated with the positive linear functional $\phi$, see
\cite{van2} and \cite{er} for details) and complete $\hat{\a0}$ in
the norm inherited from the operator-norm of $\clb(L^2(\phi))$ to
get $\hat{\cla}_r$. The definition of $\hat{\cla}$ is slightly
more complicated and involves the realization of $\hat{\a0}$ as
elements of the
 Banach $\ast$-algebra $L^1(\phi)$ (see \cite{er} and other relevant
references) and then taking the associated universal
$C^*$-completion. However, it is not important for us how the
explicit constructions of these two $C^*$-algebras are done; we
refer to \cite{van2}, \cite{er} for that; all we need is that
$\hat{\a0}$ is dense in both of them in the respective
norm-topologies.    It should also be mentioned that exactly as in
the classical case, there is a canonical surjective
$C^*$-homomorphism from $\hat{\cla}$ to $\hat{\cla}_r.$

 \newsection{Analytical assembly map  for ``proper and relatively compact" action of a discrete quantum group}

\subsection{An Analogue of proper and relatively compact action in
the quantum case}
 Let us first recall from \cite{GK}, without
proof, how one can construct an analogue of the Baum-Connes
analytic assembly map for the action of the discrete quantum group
$\cla$ on some $C^*$-algebra, under some additional assumptions on
the action, which may be called ``properness and
$\cla$-compactness", since these assumptions actually hold for   a
proper and $G$-compact action  by a discrete group $G$. Our
construction is analogous to that described in, for example,
\cite{val1},\cite{val2}, for the
  discrete group (see also \cite{meyer} for some analogous constructions in the classical context).
   We essentially translate that into our noncommutative
framework step by step, and verify that it really goes through.
However, in case $S^2$ is not identity, it is somewhat tricky to
give the correct definition of $\hat{\a0}$-valued inner product,
and prove the required properties, as one has to suitably
incorporate the modular operator $\theta$.

\vspace{5mm} Let $\clc $ be a $C^*$-algebra (possibly nonunital).
Assume furthermore that there is an action of the quantum group
$\cla$ on it, given by $\Delta_\clc :\clc \raro \clm(\clc \ot
\cla)$, which is a coassociative $C^*$-homomorphism, and assume
also that there is a dense $\ast$-subalgebra
$\clc_0$ of $\clc$ such that the following conditions are satisfied :\\
{\bf A1} $\Delta_\clc(c)(c^\prime \ot 1) \in \clc_0 \ot_{\rm alg}
\a0$ for all
$c,c^\prime \in \clc_0;$ \\
{\bf A2} $\Delta_\clc (c)(1 \ot a) \in \clc_0 \ot_{\rm alg} \a0$
for all $c
\in \clc_0,$ $a \in \a0;$\\
{\bf A3} There is a positive element $h \in \clc_0$ such that
$$(id \ot \phi)(\Delta_\clc(h^2))=1,$$ or equivalently $(id \ot
\phi)(\Delta_\clc(h^2)(c \ot 1))=c, \forall c \in \clc_0 .$

\brmrk

Assume that $X$ is a locally compact Hausdorff space equipped with
an action of a discrete group $\Gamma$ such that $X$ is
$\Gamma$-proper and $\Gamma$-compact. It is then easy to verify
that  A1,A2 and A3 hold if we take   $\clc=C_0(X),$
$\cla=C_0(\Gamma)$ and  $\clc_0= C_c(X).$ The choice of a positive
element $h$ as in A3 can be found in \cite{val1}.

\ermrk

\brmrk

 \label{ownaction}

In the quantum case, some interesting examples where A1,A2, A3 are
valid can be obtained by taking $\clc=\cla$, with the obvious
 action of $\cla$ on itself, and $\clc_0=\a0$.

 \ermrk

\subsection{The construction of analytical assembly map}

Now, our aim is to construct maps $\mu_i : KK_i^\cla(\clc,\IC)
 \raro KK_i(\IC,\hat{\cla}) \equiv K_i(\hat{\cla}),$ and $\mu_i^r :
KK_i^\cla(\clc,\IC)   \raro KK_i(\IC,\hat{\cla}_r) \equiv
K_i(\hat{\cla}_r),$ for $i=0,1$, i.e. even and odd cases. For
simplicity let us do it for $i=1$ only, the other case can be
taken care of by obvious modifications. We have chosen the
convention of \cite{val1} to treat separately odd and even cases,
instead of treating both of them on the same footing as in the
original work of Kasparov or in \cite{kk}. This is merely a matter
of notational simplicity. For the definition and properties of
equivariant KK groups $KK^\cla_.(.,.)$, we refer to  the paper by
Baaj and Skandalis (\cite{bs}) (with the easy modifications of
their definitions to treat odd and even cases separately).

Let $(U,\pi,,F)$ be a
cycle (following \cite{val1}) in $KK_1^\cla(\clc,\IC),$ i.e.\\
  (i) $U \in
\cll(\clh \ot \cla)) \cong \clm(\clb_0(\clh) \ot \cla)$ is a
unitary representation of $\cla$, where  $\clh$ is a separable
Hilbert space, i.e. $U$
is unitary and $(id \ot \Delta)(U)=U_{12}U_{13},$ $(id \ot S)(U)=U^*;$ \\
(ii) $\pi : \clc \raro \clb(\clh)$ is a nondegenerate
$\ast$-homomorphism such
that $(\pi \ot id)(\Delta_\clc(a))=U(\pi(a) \ot 1)U^*, \forall a \in \clc$;\\
(iii) $F \in \clb(\clh)$ is  self-adjoint,
$[F,\pi(c)],\pi(c)(F^2-1) \in \clb_0(\clh) \forall c \in \clc,$
and $(F \ot 1)-U(F \ot 1)U^* \in \clb_0(\clh) \ot \cla$.

\bdfn

 We say that a cycle $(U,\pi,F)$ is equivariant (or $F$ is
equivariant) if $U(F \ot 1)U^*=F \ot 1$. We say that $F$ is
properly supported if for any $c \in \clc_0,$ there are {\bf
finitely many} $c_1,...c_k, b_1,...,b_k \in \clc_0$ and
$A_1,...A_k \in \clb(\clh)$ (all depending on $c$) such that $F
\pi(c)=\sum_i \pi(c_i)A_i \pi(b_i).$

\edfn

Before we proceed further, let us make the following convention :
  for any element $A \in \clb(\clh),$ we shall denote by
$\tilde{A}$ the element $A \ot 1_{\clk}$ in $\clb(\clh \ot \clk)$.

\bthm \cite{GK}
 \label{proper} Given a cycle $(U, \pi, F)$,
 we can find  a homotopy-equivalent  cycle $(U,\pi,
F^\prime)$ such that $(U,\pi, F^\prime)$ is equivariant and
$F^\prime$ is properly supported. \ethm

{\it Sketch of Proof :}\\
Since $\pi$ is nondegenerate, we can choose a net $e_\nu$ of
elements from $\clc_0$ such that $\pi(e_\nu)$ converges to the
identity of $\clb(\clh)$ in the strict topology, i.e. in the
strong $\ast$-topology. Now, let $X_\nu:=\tilde{\pi(e_\nu)^*}U
\tilde{\pi(h)}\tilde{F}
\tilde{\pi(h)}U^*\tilde{\pi(e_\nu)}=\tilde{\pi(e_\nu)^*} (\pi \ot
id)(\Delta_\clc(h))U\tilde{F}U^* (\pi \ot
id)(\Delta_\clc(h))\tilde{\pi(e_\nu)}.$ Since by our assumption
$\tilde{e_\nu^*}\Delta_\clc(h) \in \clc_0 \ot_{\rm alg} \a0,$ and
similar thing is true for $\Delta_\clc(h)\tilde{e_\nu}$, it is
easy to see that $X_\nu$ is of the form $X_\nu=\sum_j (\pi(c_j)
\ot a_j)(U\tilde{F}U^*)(\pi(c^\prime_j) \ot a^\prime_j),$ for some
finitely many $c_j, c_j^\prime \in \clc_0$ and $a_j, a^\prime_j
\in \a0.$ Choosing a suitably large enough finite subset $I_1$ of
$I$, we can assume that all the $a_j, a^\prime_j$'s are in the
support of $e_{I_1}$, and hence it is easy to see that $X_\nu \in
\clb(\clh) \ot_{\rm alg} (e_{I_1}\a0 e_{I_1}),$ so $(id \ot
\phi)(X_\nu)$ is finite. Similarly, $(id \ot
\phi)(\tilde{\pi(e_\nu)^*}U\tilde{\pi(h^2)}U^*\tilde{\pi(e_\nu)})$
is finite, and by assumption  {\bf A3}, is equal to $(\pi(e_\nu^*
e_\nu) \ot 1).$ Now, from the operator inequality $-\|F\| 1 \leq F
\leq \| F\| 1,$ we get the operator inequality
$$ -\tilde{\pi(e_\nu)^*}U\tilde{\pi(h^2)}U^*\tilde{\pi(e_\nu)}\|F\| \leq
 X_\nu \leq \tilde{\pi(e_\nu)^*}U\tilde{\pi(h^2)}U^*\tilde{\pi(e_\nu)}\|F\|;$$
 from which it follows after applying $(id \ot \phi)$ that
$$ -\pi(e_\nu^*e_\nu) \|F \| \leq (id \ot \phi) X_\nu \leq \pi(e_\nu^*e_\nu)
\|F\|.$$
 Since $\pi(e_\nu^*e_\nu) \raro 1_{\clb(\clh)}$ in the strong operator
topology, one can easily prove by the arguments similar to those
in \cite{val1} that $(\id \ot \phi)(X_\nu)$ converges in the
strong operator topology of $\clb(\clh)$, and let us denote this
limit by $F^\prime.$ It is also easy to see that in fact
$F^\prime=(id \ot \phi)(U (\pi(h)F\pi(h) \ot 1)U^*),$ where we
have used the extended definition of $(id \ot \phi)$ on
$\clm(\clb_0(\clh) \ot \cla)$ as discussed in the previous
section.

Fix some $c \in \clc_0.$ Clearly we have $F^\prime \pi(c)=(id \ot
\phi)(U \tilde{\pi(h)}\tilde{F} \tilde{\pi(h)}U^*
\tilde{\pi(c)}).$ Now, note that
 $U \tilde{\pi(h)}\tilde{F} \tilde{\pi(h)}U^* \tilde{\pi(c)}=(\pi \ot
id)(\Delta_\clc(h))U\tilde{F}U^* (\pi \ot id)((\Delta_\clc(h)(c
\ot 1)).$ Since $\Delta_\clc(h)(c \ot 1) \in \clc_0 \ot_{\rm alg }
\a0,$ we can write it as a finite sum of the form $\sum_{ij,\al}
x^\al_{ij} \ot e^\al_{ij},$ with $x^\al_{ij} \in \clc_0,$ and
where $e^\al_{ij}$ 's are the matrix units of $\cla_\al$, as
described in the previous section, and $\al$ in the above sum
varies over some finite set $T$, say, with $i,j=1,...,n_\al.$
Thus, $  U \tilde{\pi(h)}\tilde{F} \tilde{\pi(h)}U^*
\tilde{\pi(c)} =\sum_{\al,i,j} (\pi \ot id)(\Delta_\clc(h))(1 \ot
e^\al_{ij})(Fx^\al_{ij} \ot 1).$ Since for each $\al,i,j,$
$\Delta_\clc(h)(1 \ot e^\al_{ij}) \in \clc_0 \ot_{\rm alg} \a0,$
we can write $\Delta_\clc(h)(1 \ot e^\al_{ij}) $ as a finite sum
of the form $ \sum x_p \ot a_p$ with $x_p \in \clc_0, a_p \in
\a0,$ and hence $  U \tilde{\pi(h)}\tilde{F} \tilde{\pi(h)}U^*
\tilde{\pi(c)}$ is clearly a finite sum of the form $\sum_k
\pi(c_k)A_k\pi(c^\prime_k) \ot a_k,$ with $c_k,c_k^\prime \in
\clc_0,$ $A_k \in \clb(\clh)$ and $a_k \in \a0.$ From this it
follows that $F^\prime$ is properly supported.

It is easy to show the equivariance of $F^\prime.$ Indeed,
$U(F^\prime \ot 1)U^*=(id \ot id \ot \phi)((id \ot \Delta)(U
\tilde{\pi(h)}F\tilde{\pi(h)}U^*))$ by using the fact that $(id
\ot \Delta)(U)=U_{12}U_{13}$ and $\Delta$ is a
$\ast$-homomorphism. Now, since it is easy to see using what we
have proved in the earlier section that $(id \ot id \ot \phi)((id
\ot \Delta)(X))=(id \ot \phi)(X) \ot 1,$ for $X \in
\clm(\clb_0(\clh) \ot \cla),$ from which the equivariance of
$F^\prime$ follows.

Finally, we can verify that $\pi(c)(F-F^\prime)$ is compact for $c
\in \clc_0$, hence for all $c \in \clc,$ by very similar arguments
as in \cite{val1}, adapted to our framework in a suitable way. We
omit this part of the proof, which is anyway straightforward.
\qed

\vspace{5mm}

Let $\clh_0:=\pi(\clc)\clh.$ By the fact that $F^\prime$ is
properly supported, it is clear that $F^\prime \clh_0 \subseteq
\clh_0.$  We now equip $\clh_0$ with a right $\hat{\a0}$-module
structure. Define  $$(\xi.a):=(id \ot
\psi_{\theta^{-1}S(a)\theta^{-2}})(U) \xi,$$  for $\xi \in \clh_0,
a \in \a0.$ It is useful to note that for $c \in \clb(\clh)
\ot_{\rm alg}  \a0,$ $(id \ot
\psi_{\theta^{-1}S(a)\theta^{-2}})(c)=(id \ot (\psi_{\theta a}
\circ S^{-1}))(c)$ by simple calculation using the properties of
$\psi$ and $\theta$ described in the previous section. By taking
suitable limit, it is easy to extend this for $c \in
\clm(\clb_0(\clh) \ot \cla)$, in particular for $U$. So we also
have that $\xi .a =(id \ot \psi_{\theta a} \circ S^{-1})(U) \xi.$
\bppsn \cite{GK}
 $(\xi.a).b=\xi.(a \ast b)$ for $a,b \in \a0, \xi
\in \clh_0.$ That is, $(\xi,a) \mapsto \xi .a$ is indeed a right
$\hat{\a0}$-module action.
\eppsn
{\it Proof :-}\\
  Choosing finite subsets $J,K$ of $I$ such that $\theta^{-1}S(a) \theta^{-2}
\in supp(e_K), \theta^{-1}S(b)\theta^{-2} \in supp(e_J),$ we have
that \bean
\lefteqn{(\xi.a).b}\\
&=& \sum_{\al \in J; i,j=1,..,n_\al}
U^\al_{ij}(\xi.a)\psi_{\theta^{-1}S(b)\theta^{-2}}(e^\al_{ij})\\
&=& \sum_{\al \in J;i,j=1,...,n_\al} \sum_{\beta \in
K;k,l=1,...,n_\beta} U^\al_{ij}U^\beta_{kl} \xi
\psi_{\theta^{-1}S(b)\theta^{-2}}(e^\al_{ij})\psi_{\theta^{-1}S(a)\theta^{-2}}(e^\beta_{kl})\\
&=& (id \ot \psi_{\theta^{-1}S(b)\theta^{-2}} \ot
\psi_{\theta^{-1}S(a)\theta^{-2}})(U_{12}U_{13})\xi \\
&=&  (id \ot \psi_{\theta^{-1}S(b)\theta^{-2}} \ot
\psi_{\theta^{-1}S(a)\theta^{-2}})((id \ot \Delta)(U))\xi\\
&=& (id \ot (\psi_{\theta^{-1}S(b)\theta^{-2}} \ast
\psi_{\theta^{-1}S(a)\theta^{-2}}))(U)\xi\\
&=& (id \ot \psi_{(\theta^{-1}S(b)\theta^{-2})\ast
(\theta^{-1}S(a)\theta^{-2})})(U) \xi.\\ \eean Now, by a
straightforward calculation using the properties of $\psi,$ $S$
and $\theta$ one can verify that $(\theta^{-1}S(b)\theta^{-2})\ast
(\theta^{-1}S(a)\theta^{-2})=\theta^{-1}S(a \ast b)\theta^{-2},$
which completes the proof. \qed

\vspace{5mm} For $\xi,\eta \in \clh_0,$ say of the form
$\xi=\pi(c_1)\xi^\prime,\eta=\pi(c_2)\eta^\prime,$ it is clear
that $T_{\xi \eta}(U)$ is an element of $\a0,$ since
$(\tilde{\pi(c_1^*)}U\tilde{\pi(c_2)}=(\pi \ot id)((c_1^* \ot
1)\Delta_\clc(c_2))U,$ which belongs to $  (\pi(\clc_0) \ot_{\rm
alg} \a0)\clm(\clb_0(\clh) \ot \cla) \subseteq \clb(\clh) \ot_{\rm
alg} \a0.$ We define $$ <\xi,\eta>_{\hat{\a0}}:=\theta^{-1}T_{\xi
\eta}(U) \in \hat{\a0},$$
 identifying $\a0$ as the $\ast$-algebra $\hat{\a0}$ described earlier.

We   recall here the proof in \cite{GK} of the fact  that $\clh_0$
with the above right $\hat{\a0}$-action and the $\hat{\a0}$-valued
bilinear form $<.,.>_{\hat{\a0}}$ is indeed a pre-Hilbert
$\hat{\a0}$-module.

Define $\Sigma : \clh_0 \raro \clf_0$ by $$ \Sigma(\xi):= ((\pi(h)
\ot \theta^{-1})U )\xi,$$
 for $\xi \in \clh_0.$ Note that by writing $\xi=\pi(c)\xi^\prime$ for some
$\xi^\prime \in \clh, c \in \clc_0,$ we have that $((\pi(h) \ot
1)U) \xi=
 ((\pi \ot id)((h \ot 1)\Delta_\clc(c))U)\xi^\prime,$ and since $   (\pi \ot
id)((h \ot 1)\Delta_\clc(c))U \in \pi(\clc_0) \ot_{\rm alg} \a0,$
the range of $\Sigma$ is clearly in $\clh \ot_{\rm alg} \a0.$ It
follows  that $\Sigma$ is in fact a module map and preserves the
bilinear form $<.,.>_{\hat{\a0}}$ on $\clh_0.$

\bppsn
 \cite{GK}
\label{sigma} For $\xi,\eta \in \clh_0, a \in \a0,$ we have that\\
(i) $\Sigma(\xi .a)=\Sigma(\xi)a.$\\ (ii)
$<\Sigma(\xi),\Sigma(\eta)>=<\xi,\eta>_{\hat{\a0}}.$
 \eppsn
{\it Proof :-}\\
(i) Choose suitable finite set indexed by $p$ such that
$(\tilde{\pi(h)}U)\xi=\sum_p \pi(h)U_1^{(p)} \xi \ot U^{(p)}_2,$
where $U^{(p)}_1 \in \clb(\clh), U^{(p)}_2 \in \a0,$ and  also
$\sum_p U_1^{(p)^*}  \ot U^{(p)^*}_2=\sum_p U_1^{(p)}  \ot
S(U^{(p)}_2).$ Using the facts that $\Delta(\theta)=\theta \ot
\theta,$ $S^{-1}(\theta)=\theta^{-1}$ and that
$\psi(b\theta)=\psi(\theta b) \forall b \in \a0,$ and also the
easily verifiable relation $\psi_a \circ
S^{-1}=\psi_{\theta^{-1}S(a)\theta^{-1}}$ for $a \in \a0,$ we have
that \bean
\lefteqn{ \Sigma(\xi)a}\\
&=&   \sum_p \pi(h)U_1^{(p)} \ot (id \ot
\psi_{\theta^{-1}S(a)\theta^{-1}})(\Delta(\theta^{-1}U_2^{(p)}))\\
&=& \sum_p \pi(h)U_1^{(p)} \ot (id \ot
\psi_{\theta^{-1}S(a)\theta^{-1}})( (\theta^{-1} \ot
\theta^{-1})\Delta(U_2^{(p)}))\\
&=& (\pi(h) \ot \theta^{-1} \ot
\psi_{\theta^{-1}S(a)\theta^{-2}})((\id \ot
\Delta)(U\xi))\\
&=& (\pi(h) \ot \theta^{-1} \ot
\psi_{\theta^{-1}S(a)\theta^{-2}})((id \ot
\Delta)(U) \xi)\\
&=& (\pi(h) \ot \theta^{-1} \ot
\psi_{\theta^{-1}S(a)\theta^{-2}})(U_{12}U_{13})\xi \\
&=& (\pi(h) \ot \theta^{-1})(U) ((id \ot
\psi_{\theta^{-1}S(a)\theta^{-2}})(U)\xi)\\
&=& (\pi(h) \ot \theta^{-1})(U)(\xi.a)\\
&=& \Sigma(\xi a).
\eean

(ii) Choosing suitable finite index sets
as explained before, such that
 $(\pi(h) \ot 1)U \xi=\sum_p U_1^{(p)} \ot U_2^{(p)},$ with $\sum_p
U_1^{(p)} \ot S(U_2^{(p)})=\sum_p U_1^{(p)^*} \ot U_2^{(p)^*},$
and similarly for $(\pi(h) \ot 1)U \eta$ with the index $p$
replaced by say $q$, we can write \bean
\lefteqn{<\Sigma(\xi),\Sigma(\eta)>}\\
&=& \sum_{p,q} <U_1^{(p)}\xi,
\pi(h^2)U_1^{(q)}\eta>(\theta^{-1}U_2^{(p)})^\sharp \ast
(\theta^{-1}U_2^{(q)})\\
&=& \sum_{p,q} <\xi, U_1^{(p)^*}
\pi(h^2)U_1^{(q)}\eta>(\theta^{-2}S^{-1}(\theta^{-1})S^{-1}(U_2^{(p)^*}))
\ast (\theta^{-1}U_2^{(q)})\\
&=& \sum_{p,q} <\xi, U_1^{(p)^*}
\pi(h^2)U_1^{(q)}\eta>(\theta^{-1}S^{-1}(U_2^{(p)^*}))
\ast (\theta^{-1}U_2^{(q)})\\
&=& \sum_{p,q} <\xi,
U_1^{(p)}\pi(h^2)U_1^{(q)}\eta>\theta^{-1}(U_2^{(p)}\ast
 U_2^{(q)})\\
&=& \sum_{p,q} <\xi, U_1^{(p)}\pi(h^2)U_1^{(q)}\eta> (\phi \ot
\theta^{-1})((U_2^{(p)}\ot 1)(S \ot id)(\Delta(U_2^{(q)})))...(1),
 \eean
using the fact that $\sum_p U_1^{(p)^*} \ot U_2^{(p)^*}=\sum_p
U_1^{(p)} \ot S(U_2^{(p)})$ and the simple   observation  that
$(\theta^{-1}x) \ast (\theta^{-1}y)=\theta^{-1}(x \ast y).$ Now,
\bean
\lefteqn{\sum_q \pi(h^2)U_1^{(q)}\eta \ot ((S \ot id)(\Delta(U_2^{(q)})))}\\
&=& (\pi(h^2) \ot S \ot id)((id \ot \Delta)(U\eta))\\
&=& (\pi(h^2) \ot S \ot id)((id \ot \Delta)(U) \eta)\\
&=& (\pi(h^2) \ot id \ot id)((U^*)_{12}U_{13})\eta ....(2). \eean
 Thus, from (1) and (2), $<\Sigma(\xi),\Sigma(\eta)>=
(T_{\xi \eta} \ot \phi \ot \theta^{-1})((U(\pi(h^2) \ot 1)U^* \ot
1)U_{13})=\theta^{-1}T_{\xi \eta}(U),$ since $(id \ot \phi)(U
\pi(h^2) U^*)=1.$ This completes the proof.
\qed

Note that from the above proposition it follows in particular that
$<\xi,\eta a>_{\hat{\a0}}=<\Sigma(\xi),\Sigma(\eta
a)>=<\Sigma(\xi),\Sigma(\eta)a>=<\Sigma(\xi),\Sigma(\eta)>\ast a
=<\xi,\eta>_{\hat{\a0}} \ast a.$ Similarly,
$<\xi,\eta>_{\hat{\a0}}^\sharp=<\eta,\xi>_{\hat{\a0}},$ and
$<\xi,\xi>$ is a nonnegative element in the $\ast$-algebra
$\hat{\a0}$, since $<.,.>$ on $\clf_0$ is a nonnegative definite
form.

 Given any $C^*$-algebra which contains $\a0$ as a dense $\ast$-subalgebra, we
can complete $\clf_0$ w.r.t. the corresponding norm to get a
Hilbert $C^*$-module in which $\clf_0$ sits as a dense submodule.
Let us denote by $\clf$ and $\clf_r$ the Hilbert $\hat{\cla}$ and
$\hat{\cla}_r$-modules respectively obtained in the above
mentioned procedure, by considering $\a0$ as dense
$\ast$-subalgebra of $\hat{\cla}$ and $\hat{\cla}_r$ respectively.
The corresponding completions of $\clh_0$ will be denoted by
$\cle$ and $\cle_r$ respectively. By construction, $\Sigma$
extends to an isometry from $\cle$ to $\clf$ and also from
$\cle_r$ to $\clf_r$. We denote both these extensions by the same
notation $\Sigma$, as long as no confusion arises. Clearly, $\cle
\cong \Sigma \cle \subseteq \clf$ as closed submodule, and similar
statement will be true for $\cle_r$ and $\clf_r.$

 \bthm
 \cite{GK}
\label{compact}
Let $T \in \clb(\clh)$ be equivariant, i.e. $U(T
\ot 1)U^*=T \ot 1,$ and also assume that it satisfies the
following condition
 which is slightly weaker than being properly supported :\\
For $c \in \clc_0,$ one can find $c_1,...,c_m \in \clc_0,
A_1,...,A_m \in \clb(\clh)$ (for some integer $m$) such that $T
\pi(c)=\sum_k
\pi(c_k)A_k.$\\
  Then we have the following :\\
(i) $T(\xi a)=(T\xi)a$ $ \forall a \in \a0,$ and thus $T$ is a
module map on the $\hat{\a0}$-module $\clh_0.$ Furthermore, if $T$
is self-adjoint in the sense of Hilbert space, then
$<\xi,T\eta>_{\hat{\a0}}=<T\xi,\eta>_{\hat{\a0}}$ for $\xi,\eta
\in \clh_0.$ \\ (ii) $T$ is continuous in the norms of  $\cle$ as
well as $\cle_r$, thus admits continuous extensions on both $\cle$
and $\cle_r$. We shall denote these extensions by $\clt$ and
$\clt_r$ respectively.\\ (iii) If $T \pi(h)$ is compact in the
Hilbert space sense, i.e. in $\clb_0(\clh),$ then $\clt$ and
$\clt_r$ are compact in the Hilbert module sense.
\ethm
{\it Proof :}\\
(i) is obvious from the defintion of the right $\hat{\a0}$ action,
the definition of $<.,.>_{\hat{\a0}}$, and the equivariance of
$T$. Let us prove (ii) and (iii) only for $\clt$, as the proof for
$\clt_r$ will be exactly the same. In fact, it is enough to show
that $\Sigma \clt \Sigma^*$ is continuous on $\clf$, and is
compact if $T \pi(h)$ is compact in the Hilbert space sense. Let
us introduce the following notation : for $X \in \clm(\clb_0(\clh)
\ot \cla), a \in \a0, \eta \in \clh,$ define $X \ast b :=(id \ot
id \ot \psi_b \circ S^{-1})((id \ot \Delta)(X)),$ and $X \ast (
\eta \ot a):=(X \ast a)\eta.$ Note that clearly $X \ast a \in
\clm(\clb_0(\clh) \ot \cla)$, so $(X \ast a) \eta$ makes sense.
Now, we observe using the equivariance of $T$ and the explicit
formula for $\Sigma^*$ derived earlier that for $\eta \in \clh, a
\in \a0,$
\bean \lefteqn{\Sigma \clt \Sigma^*(\eta \ot a)}\\
&=& (\pi(h) \ot 1)U \beta, \eean where $\beta \in \clh$ is given
by $\beta=(id \ot
\psi_{\theta^{-1}S(a)\theta^{-2}})(U)(T\pi(h)\eta).$  Now, by
using the fact that $(id \ot \Delta)(U)=U_{12}U_{13},$ it follows
by a straightforward computation that \bean
\lefteqn{ (1 \ot \theta^{-1})U \beta}\\
&=& (id \ot \theta^{-1} \ot (\psi_{\theta a} \circ S^{-1}))((id
\ot \Delta)(U) (T \pi(h) \eta). \eean
   But $ \psi_{\theta a}( S^{-1}(b))=\psi(\theta a
S^{-1}(b))=\psi(aS^{-1}(b)\theta)=\psi(aS^{-1}(\theta^{-1}b))=(\psi_a
\circ S^{-1})(\theta^{-1}b),$ and hence $(id \ot \theta^{-1} \ot
(\psi_{\theta a} \circ S^{-1} ))((id \ot \Delta)(U))=(id \ot id
\ot (\psi_a \circ S^{-1}))((id \ot \Delta)((1 \ot
\theta^{-1})U))=((1 \ot \theta^{-1})U) \ast a.$ From this, it is
clear that
$$ (\Sigma \clt \Sigma^*)(\eta \ot a)=((\pi(h) \ot \theta^{-1})U(T \pi(h) \ot
1))\ast (\eta \ot a).$$
 Now, note that $T \pi(h)=\sum_{k=1}^m
\pi(c_k)A_k,$ for some $c_1,...,c_m \in \clc_0, A_1,...,A_m \in
\clb(\clh),$ and so we have  $(\pi(h) \ot \theta^{-1})U(T \pi(h)
\ot 1) =\sum_k (1 \ot \theta^{-1})(\pi \ot id)((h \ot
1)\Delta_\clc(c_k))U(A_k  \ot 1).$  But $(h \ot
1)\Delta_\clc(c_k)$ is in $\clc_0 \ot_{\rm alg} \a0$ for each
$k=1,..m,$  and thus   $(\pi(h) \ot \theta^{-1})U(T \pi(h) \ot 1)
\in \clb(\clh) \ot_{\rm alg} \a0$ clearly. Choosing some large
enough
 finite subset  $J$ of $I$ such that $(\pi(h) \ot \theta^{-1})U(T \pi(h) \ot
1) =\sum_{\al \in J, ij=1,..,n_\al} B^\al_{ij} \ot e^\al_{ij}$,
(with $B^\al_{ij} \in \clb(\clh)$), it is clear that $\Sigma \clt
\Sigma^*=\sum_{\al \in J, ij=1,...,n_\al} B^\al_{ij} \ot
L_{e^\al_{ij}}$, where for $x \in \a0,$ $L_x : \a0 \raro \a0$ with
$L_x(a)=x \ast a.$ As $L_x$ is a norm-continuous map on
$\hat{\cla},$ the above finite sum shows that $\Sigma \clt
\Sigma^*$ indeed admits a continuous extension on the Hilbert
$\hat{\cla}$-module $\clf$. This proves (ii).

Furthermore, since $\clk(\clh \ot \hat{\cla}) \cong \clb_0(\clh)
\ot \hat{\cla}$, where $\clk(E)$ means the set of compact (in the
Hilbert module sense) opeartors on the Hilbert module $E$, it is
easy to see that $\Sigma \clt \Sigma^*$ is compact on $\clf$ if
$B^\al_{ij}$'s are compact on the Hilbert space $\clh.$ Now,
$B^\al_{ij}=(id \ot \phi^\al_{ij}) ((\pi(h) \ot \theta^{-1})U(T
\pi(h) \ot 1))=\pi(h).(id \ot \phi^\al_{ij})((1 \ot
\theta^{-1})U)T \pi(h),$ where $\phi^\al_{ij} $ is the functional
on $\a0$ which is $0$ on all $e^\beta_{kl}$ except $\beta=\al,
(kl)=(ij),$  with $\phi^\al_{ij}(e^\al_{ij})=1.$ It follows that
$B^\al_{ij}$ 's are all compact if $T \pi(h)$ is so, which
completes the proof.
\qed

Now, let us come to the construction of the Baum-Connes maps
$\mu_1 : KK_1^\cla(\clc,\IC) \raro KK_1(\IC,\hat{\cla})$ and
$\mu_1^r : KK_1^\cla(\clc,\IC) \raro KK_1(\IC,\hat{\cla}_r).$ Let
us do it only for $\mu_1$, as the case of $\mu^r_1$ is similar,
and in fact $\mu^r_1$ will be the compositon of $\mu_1$ and the
canonical map from $KK_1(\IC , \hat{\cla})$ to
$KK_1(\IC,\hat{\cla}_r)$ induced by the canonical surjective
$C^*$-homomorphism from $\hat{\cla}$ to $\hat{\cla}_r.$ Note that
an element of $KK_1(\IC,\hat{\cla}) \cong K_1(\hat{\cla})$  is
given by the suitable homotopy class $[E,L]$ of a pair of the form
$(E, L)$, where $E$ is a Hilbert $\hat{\cla}$-module and $L \in
\cll(E)$ (the set of adjointable $\hat{\cla}$-linear maps on $E$)
such that $L^*=L,$ $L^2-1$ is compact in the sense of Hilbert
module. For more details, see for example \cite{kk}.

\bthm
 \cite{GK}
\label{defn} Given a cycle $(U, \pi, F) \in KK_1^\cla(\clc,\IC),$
let $F^\prime \equiv F^\prime_h$ be the equivariant and properly
supported operator as constructed in \ref{proper}, with a given
choice of $h$ as in that theorem. Then the continuous extension of
$F^\prime_h$ on the Hilbert module $\cle$ (as described by  the
Theorem \ref{compact}), to be denoted  by say $\clf^\prime_h$,
satisfies the conditions that $(\clf^\prime_h)^*=\clf^\prime_h$
(as module map), and  $(\clf^\prime_h)^2-I$ is compact on $\cle$.
Define  $$ \mu_1((U,\pi,F)):=[\cle, \clf^\prime_h] \in
KK_1(\IC,\hat{\cla}) \cong K_1(\hat{\cla}).$$ In fact,
$[\cle,\clf^\prime_h]$  is independent (upto operatorial homotopy)
of the choice of $h$.
\ethm
{\it Proof :-}\\
  Since $F^\prime_h$ is equivariant and properly supported, it is clear that
$T_h:=(F^\prime_h)^2-1$ is equivariant and for any $c \in \clc_0,$
there are finitely many $c_1,...,c_m \in \clc_0, A_1,...,A_m \in
\clb(\clh)$ such that $T_h \pi(c)=\sum_k \pi(c_k)A_k.$
Furthermore, by the Theorem \ref{proper}, we have that
$\pi(c)T_h$, and hence $T_h \pi(c)$ is compact operator on $\clh$
 for every $c \in \clc.$ So, in particular, $T_h \pi(h)$ is compact. By
Theorem \ref{compact}, it follows that the continuous extension of
$T_h$ on $\cle$ is compact in the sense of Hilbert modules.
Furthermore, the fact that $(\clf^\prime_h)^*=\clf^\prime_h$ is
clear from (i) of the Theorem \ref{compact}. So,
$[\cle,\clf^\prime_h] \in KK_1(\IC,\hat{\cla}).$ Furthermore, as
we can see from the proof of   the Theorem \ref{proper},
$(F^\prime_h-F)\pi(c) \in \clb_0(\clh)$ $\forall c \in \clc_0,$
and so for $h,h^\prime$ satisfying {\bf A3}, we have
$(F^\prime_h-F^\prime_{h^\prime})\pi(c) \in \clb_0(\clh),$ and
hence by Theorem \ref{compact},
$\clf^\prime_h-\clf^\prime_{h^\prime}$ is compact in the Hilbert
module sense. Thus, for each $t \in [0,1],$  setting
$\clf(t):=t\clf^\prime_{h^\prime}+(1-t)\clf^\prime_h,$ we have
that $\clf(t)^2-I $ is compact on $\cle$, and this gives a
homotopy in $KK_1(\IC,\hat{\cla})$ between $[\cle,\clf^\prime_h]$
and $[\cle,\clf^\prime_{h^\prime}].$
\qed

\subsection{Examples and computations}
{\large {\bf Example 1 : Finite Dimensional Quantum Groups}}\\

Assume that $\cla$ is a finite dimensional quantum group (hence
both discrete and compact). Note that in this case $\hat{\cla}$
coincides with $\hat{\cla_r}$, $\phi=\psi$ is a tracial state,
$S^2=id,$ $h=1$ and $\theta=1$. We would now like to prove the
following

\bppsn

\label{finite}

The analytical assembly map gives an isomorphism between the
abelian groups $KK^\cla_i(\IC,\IC)$ and $KK_i(\IC,\hat{\cla})$ for
$i=0,1.$

\eppsn

{\it Proof :-}\\
 It is enough to do the case $i=0,$ since for $i=1,$ both
$KK^\cla_1(\IC,\IC)$ and $KK_1(\IC,\hat{\cla})$ are easily seen to
be $\{0\}.$ Let us assume that the set of irreducible
representation of the quantum group $\cla$ has cardinality $N$,
and let us index them by say $\sigma_1,...\sigma_N.$ Then,
$\hat{\cla} \cong \oplus_{i=1}^N M_{d_i},$ where $d_i$ denotes the
dimension of the representation $\sigma_i.$ Consider a cycle $[U,
\clh, F, \gamma]$ in $KK_0^\cla(\IC, \IC), $ where
 $\gamma$ is the
  grading operator on $\clh,$ and without loss of generality assume that $F$ and $\gamma$ are equivariant,
   $F$ is self-adjoint and $F^2-1 \in \clb_0(\clh).$ We can decompose (by applying the
Peter-Weyl Theorem for  compact quantum groups to  $\cla$) $\clh$
as $\clh=\oplus_i \IC^{d_i} \ot \clh_i,$ so that $U=\oplus_i
U^{(\sigma_i)} \ot I_{\clh_i},$ $F=\oplus_i I \ot F_i,$
$\gamma=\oplus_i I \ot \gamma_i$, where $U^{(\sigma)} \equiv ((
U^{\sigma}_{kl} ))_{k,l=1}^{d_\sigma}$ denotes the $d_\sigma
\times d_\sigma$ unitary matrix corresponding to the irreducible
representation of type $\sigma$ and $\clh_i$ is separable Hilbert
space, $F_i \in \clb(\clh_i)$ is a self-adjoint operator acting on
$\clh_i$ such that $F_i^2-1_{\clh_i} \in \clb_0(\clh_i),$
$\gamma_i$ is a grading operator acting on $\clh_i$ with $F_i
\gamma_i =- \gamma_i F_i.$ In other words, $\beta_i \equiv
(\clh_i, F_i, \gamma_i)$ can be considered as a cycle in
$KK_0(\IC, \IC)=K_0(\IC) \cong \IZ,$ and it is clear that the map
$[U,\clh, F, \gamma] \mapsto (\beta_1,...,\beta_N) \in \oplus _i
KK_0(\IC, \IC) \cong KK_0 (\hat{\cla}) \cong \IZ^N$ is an
isomorphism.
 We claim that
 the analytical assembly map $\mu^r_0=\mu_0$ is nothing but the above
  association.

 Recall from earlier sections that as a $C^*$-algebra $\hat{\cla}$ is nothing but $\cla$ (here the dimension is finite, so no further completion is needed)
  with   the product given by $(a,b) \mapsto a \ast b.$ It is not difficult to verify that
   $U^{(\sigma_i)}_{kl} \ast U^{(\sigma_j)}_{pq}=\delta_{ij} \delta_{lp} U^{(\sigma_i)}_{kq},$ upto some constant
    multiple,  (as finite dimensional quantum groups
    are unimodular, the matrix elements $U^{(\sigma)}_{kl}$ are orthogonal w.r.t. the inner product coming from the Haar
     state) and
    we can identify the $\ast$-algebra  generated by $U^{(\sigma_i)}_{kl}$ with
     the component $M_{d_i}$ in $\hat{\cla}=\oplus_j M_{d_j},$   identifying  (upto some constant)
       $U^{(\sigma_i)}_{kl}$ with $|e^i_k><e^i_l|,$ where $\{ e^i_k ,k=1,2,... \}$ is an orthonormal basis
       of $\IC^{d_i}.$   Let us now
  look at the action of
  $\mu_0$ on the cycle
$[ U^{(\sigma_i)} \ot I, \IC^{d_i} \ot \clh_i, I \ot F_i, I \ot
\gamma_i].$  For $\xi, \eta \in \clh_i, $ and $k,l \in \{
1,2,...,d_i \},$ it is easy to see that $< e^i_k \ot \xi, e^i_l
\ot \eta>_{\hat{\cla}}=<\xi, \eta>|e^i_k><e^i_l| $, and  the
$\hat{\cla}$-module obtained by the map $\mu_0$ is $E_i \ot
\clh_i,$ where $E_i$ denotes the $\IC^{d_i}$ with the canonical
$M_{d_i}$-module given by $\xi.a:=a^\prime \xi,$ $<\xi, \eta>:=
\sum_{kl} \bar{\xi_k} \eta_l |e^i_k><e^i_l|,$ for $\xi=\sum_k
\xi_k e^i_k, \eta=\sum_l \eta_l e^i_l \in \IC^{d_i}.$  Here
$a^\prime$ denotes the transpose of the matrix $a.$ Clearly, $E_i$
can be viewed as $\hat{\cla}$-module with the same inner product
and the right action extended by $\xi.b=0$ for $b $ in any
component other than $M_{d_i}.$ Thus, $\mu_0$ takes $[
U^{(\sigma_i)} \ot I, \IC^{d_i} \ot \clh_i, I \ot F_i, I \ot
\gamma_i]$ to $(E_i \ot \clh_i, I \ot F_i, I \ot \gamma_i)$, which
can be identified with the cycle $(\clh_i, F_i, \gamma_i)$ of
$KK_0(\IC, M_{d_i})=KK_0(\IC, \IC).$ This completes the proof.
 \qed

{\large {\bf Example 2 :-}}

Let us consider a general (possibly infinite dimensional)
separable discrete quantum group $\cla$. For simplicity, assume
that $\hat{\cla}_r$ coincides with $\hat{\cla}$ (there are many
interesting quantum groups satisfying this, e.g. the dual of
$SU_q(n)$, $n=2,3,...,$ $q$ positive). Consider $\clc=\cla, $ with
its canonical action on itself, denoted  by $\Delta.$ Take
 $\clc_0=\a0$, as mentioned in the Remark \ref{ownaction}.
   Let us recall the notation and discussion on general discrete
quantum groups in Section 3. The dual $\hat{\cla}$ is a compact
quantum group, and from the general theory it is known that the
GNS space for $(\cla, \phi)$ is isomorphic as a Hilbert space with
the GNS space of the Haar state on $\hat{\cla}$ via a
quantum-analogue of the classical Fourier transform. Let $\clh$ be
the GNS space of $(\cla, \phi).$ It is known that both $\cla$ and
$\hat{\cla}$ can be represented as subalgebras of $\clb(\clh)$.
Furthermore, there is a unitary $U$ such that $\Delta(a)=U (a \ot
I) U^*$ for $a \in \cla,$ and in fact this $U$ is given by $U(a
\ot b)=\Delta(a)(b \ot 1)$ for $a,b \in \a0.$ It can be verified
by direct computation, following step by step the construction of
$\mu_i,i=0,1$ discussed in the previous subsection that the
Hilbert $\hat{\cla}$-module arising out of this construction can
be identified with the trivial module $\hat{\cla}$ itself, with
its own right action and the canonical inner product $<x,y>=x^*y.$
 Let us see this for the case when $S^2=id,$ just for sake of
 computational simplicity, keeping in mind that the general case
 can be treated essentially in a similar manner. Note that in this case
  $\theta=1, \psi=\phi,$ and we have  for
 $\xi  \in \a0,$ viewed as an element in $L^2(\cla,\phi), $ $a \in
 \a0,$
  $\xi.a=(id \ot \psi_a \circ S^{-1})(U \xi)=(id \ot \phi_a \circ
  S )(\Delta(\xi))=\xi \ast a.$ Similarly, $<\xi, \eta>=\xi^\sharp
   \ast \eta$ can be easily seen. This proves our claim.

        As a particular example of the above, we can consider
        $\cla=\hat{SU}_q(2)$, $q$ positive. Fix a generator
        of $K_0(\IC)=\IZ $, realized as $(\clh_0, F, \gamma)$, say,
        and then consider $\clh=L^2(\cla, \phi) \ot \clh_0,$ and a
        cycle $[U \ot I, \clh, I \ot F, I \ot \gamma] \in
        KK_0^\cla(\cla, \IC).$ Here $\cla$ is represented as $a
        \mapsto (a \ot I_{\clh_0})$ in $\clb(\clh)$.  The
        analytical assembly map $\mu_0$ will take this cycle to
        the canonical generator of $K_0(SU_q(2))=\IZ,$ since the
        inclusion of $1 \in SU_q(2)$ is an isomorphism of
        $K_0$-groups.

\newsection{Formulation of the conjecture}

\subsection{Formulation}
Now we shall formulate the Baum-Connes' conjecture for $\cla$,
following the approach (for a classical discrete group) taken by
Cuntz in \cite{cuntz}. We shall use the notation introduced and
used in Section 2.

Fix any finite subset $F$ of $I$. Define a $C^*$-algebra $\cle_F$
as follows (note that for $n=1,2,..$ we denote by $\Delta^n
:\clm(\cla) \raro \clm(\cla^{\ot^{n+1}})$ the map $...(\Delta \ot
id \ot id)(\Delta \ot id) \Delta.$) :

\bdfn
  \label{univc*}
  Let $\cle_F$ be the universal $C^*$-algebra generated by
  elements $\lambda(a),a \in \cla_0$ satisfying the following :\\
  (i) $\cla_0 \ni a \mapsto \lambda(a) \in \cle_F$ is completely
  positive,\\
  (ii) $(\lambda \ot id)(\Delta^n(a)) (\lambda \ot id)(\Delta^n(x_1))...(\lambda \ot id)(\Delta^n(x_k))
  (\lambda \ot id)(\Delta^n(b))=0 \forall x_1,...,x_k \in \cla_0$ whenever
 $a$ and $b$ are such that $(1 \ot e_F)(\sigma \circ \Delta(a))(b
 \ot 1)=0
 ,$ (where $\sigma$ is the flip automorphism on $\clm(\cla \ot \cla)$),\\
 (iii) $\sum_\al \lambda(e_\al)\lambda(a)=\lambda(a)$ for all $a
 \in \cla_0.$ \\
 Note that the sum in (iii) above is actually a finite sum, as a consequence of
 (ii).

 \edfn

 We shall define a coassociative coaction $\Delta_F : \cle_F \raro
 \clm(\cle_F \ot \cla)$. We need the following result for doing
 this.
 \bppsn
 \label{coaction}
 Define $\tilde{\lambda} : \cla_0 \raro \clm(\cle_F \ot \cla)$ by
 $\tilde{\lambda}=(\lambda \ot id) \circ   \Delta.$  Then (i)-(iii)
 in the definition of $\cle_F$ are valid if we replace $\lambda$
 by $\tilde{\lambda}$; so by the universal property of $\cle_F$,
 there is a well-defined $C^*$-homomorphism $\Delta_F : \cle_F
 \raro \clm(\cle_F \ot \cla)$ satisfying
 $\Delta_F(\lambda(a))=\tilde{\lambda}(a)$ for all $a \in \cla_0.$
 \eppsn
 {\it Proof :-}\\
 Since $\lambda$ is clearly a nondegenerate, bounded completely positive
 ( to be abbreviated as CP)
 map (nondegeneracy  follows from (iii) of the definition, which says that $\sum_\al \lambda(e_\al)=1$ in
  the strict topology of $\clm(\cle_F)$),
 $(\lambda \ot id)$ makes sense on $\clm(\cla \ot \cla)$. Clearly,
 $\tilde{\lambda}$ is completely positive on $\cla_0.$ Now,
  choose $a,b \in \cla_0$ such that $(1 \ot e_F)(\sigma \circ
  \Delta(a))
  (b \ot 1)=0.$ We have,
  $(\tilde{\lambda} \ot id)(\Delta^n(a))(\tilde{\lambda}\ot id)(\Delta^n(x_1))...
  (\tilde{\lambda} \ot id)(\Delta^n(b))=
   (\lambda \ot id)(\Delta^{n+1}(a))(\lambda\ot id)(\Delta^{n+1}(x_1))...
  (\lambda \ot id)(\Delta^{n+1}(b))=
     0$ for
  all $x_1,...,x_k \in \cla_0.$ Finally, it is easy to see that
  the condition (iii) holds for $\lambda$ replaced by
  $\tilde{\lambda}.$
  \qed

  It is easy to show that $\Delta_F$ is indeed coassociative, since $\Delta$ is coassociative
   and  we have, $(\Delta_F \ot id)(\lambda \ot id)=(\lambda \ot id \ot id)(\Delta \ot id)$,
    from which it is immediate that $(\Delta_F \ot
    id)(\Delta_F(\lambda(a)))=(id \ot
    \Delta_F)(\Delta_F(\lambda(a)))$. Now we shall show that
     the $\cla$-coaction on $\cle_F$ satisfies the conditions
     A1,A2, A3 mentioned earlier by us (see 3.1), if we choose large enough $F$ such that it contains
      $\al_0$ mentioned in the earlier section.

     \bthm
     A1,A2,A3 are valid for the coaction $\Delta_F : \cle_F \raro
     \clm(\cle_F \ot \cla)$.
     \ethm
     {\it Proof :-}\\
     Let us take $\clc=\cle_F$ and let $\clc_0$ be the
     $\ast$-algebra (no completion) generated by $\lambda(a), a
     \in \cla_0.$ Take $a,b \in \cla_0,$ and let $c=\lambda(a) \in
     \clc_0.$ From the general theory of discrete quantum groups,
     we know that $\Delta(a)(1 \ot b) \in \cla_0 \ot_{\rm alg}
     \cla_0.$ Thus, $\Delta_F(\lambda(a))(1 \ot b)=(\lambda \ot
     id)(\Delta(a))(1 \ot b)=(\lambda \ot id)((\Delta(a)(1 \ot b))
      \in (\lambda \ot id)(\cla_0 \ot_{\rm alg} \cla_0)=\clc_0 \ot_{\rm alg}
      \cla_0.$ This proves A2. To prove A1, we note from
      \cite{van} that given $\beta,\gamma \in I,$ there are only finite many $\al \in I$ such
      that $\Delta(e_\al)(e_\beta \ot e_\gamma)\ne 0.$ Thus,
      depending on $b$ and $F$, we can choose a finite subset $J$
      of $I$ such that $\forall \al $ not in $J$, one has $(1 \ot
      e_F)(\Delta(e_\al)(b \ot 1))=0.$ Hence for any $x \in
      \cla_0,$ we have that $(1 \ot e_F)(\Delta(xe_\al)(b \ot
      1))=\Delta(x)\Delta(e_\al)(b \ot e_F)=0,$ for any $\al$ not
      belonging to $J$ (as $e_F$ is a central element in
      $\clm(\cla \ot \cla).$) From this, it follows that
      $\lambda(xe_\al)\lambda(b)=0$ for $\al$ not in $J$, and for
      any $x \in \cla_0,$ and hence for any $x \in \clm(\cla)$, since  $\lambda$,
      being a nondegenerate CP
       map on $\cla_0$, has a strictly continuous extension to
       $\clm(\cla)$. Thus, $ (\lambda \ot id)(X)(\lambda(b) \ot 1)=
       (\lambda \ot id)(X(e_J \ot 1))(\lambda(b) \ot 1)$ for all
       $X \in \clm(\cla \ot \cla).$ In particular, for $a,b \in
       \cla_0,$
       $(\lambda \ot id)(\Delta(a))(\lambda(b)
       \ot 1)=(\lambda \ot id)(\Delta(a)(e_J \ot 1))(\lambda(b)
       \ot 1) \in \clc_0 \ot_{\rm alg} \cla_0$, since
       $\Delta(a)(e_J \ot 1) \in \cla_0 \ot_{\rm alg} \cla_0.$
       This completes the proof of A1. Finally, by taking
       $h=\lambda(h_0),$ where $h_0$ is as in the earlier section,
        we see that A3 is satisfied, as $(id \ot \phi)(\Delta(h_0))=\phi(h_0)1=1$.
\vspace{2mm}

  Applying the results of \cite{GK}, as recalled already, we have
   canonical homomorphisms
 $\mu^F_i :
KK_i^\cla(\cle_F,\IC) \raro KK_i(\IC,\hat{\cla})$ and $\mu_i^{r,F}
: KK_i^\cla(\cle_F,\IC) \raro KK_i(\IC,\hat{\cla}_r),$ $i=0,1,$
and for every $F$ as above. Then, taking the inverse limit w.r.t.
 a  family of $F$ increasing to $I$, we get maps $\mu_i :
lim_F KK_i^\cla(\cle_F,\IC) \raro KK_i(\IC,\hat{\cla})$ and
$\mu_i^r :lim_F KK_i^\cla(\clc,\IC) \raro KK_i(\IC,\hat{\cla}_r).$

\vspace{3mm}
 {\Large {\bf 5.1.4.  conjecture :}}  $\mu_i^r$ is an isomorphism for $i=0,1.$
\vspace{3mm}
\subsection{Consistency with the classical formulation and
verification for some quantum groups}

 We first point out that in the classical case, this is
equivalent to Baum-Connes' conjecture, using the results in
\cite{cuntz}.

\bppsn

If we take $\cla=C_0(G)$, where $G$ is a discrete group, then the
conjecture 5.1.4 above is equivalent to the classical Baum-Connes
Conjecture for the group $G$.

\eppsn

{\it Proof :-}\\
 Suppose that   $F$ is some finite
subset of $G$ containing the identity element $e$, say, and
$\cla=C_0(G), \cla_0=C_c(G).$ Since $\cla$ is commutative,
complete positivity coincides with positivity, and hence any CP
map $\lambda$ from $\cla$ to some other $C^*$ algebra $\clb$ is
determined by a map from $G$ to the positive cone of $\clb$ given
by $g \mapsto \lambda_g \equiv \lambda(\delta_g) \ge 0$, where
$\delta_g$ is the indicator function at the point $g \in G.$ Thus,
in this case, $\lambda(f)=\sum_g f(g) \lambda_g$ for $f \in
C_c(G).$ Furthermore, the condition $(1 \ot e_F)(\sigma \circ
\Delta(f))(\psi \ot 1)=0$ for some pair $(f,\psi)$ is equivalent
to $f(hg)\psi(g)=0 \forall h \in F,$ i.e. $f(s) \psi(t)=0$
whenever $st^{-1} \in F.$ In particular, given two fixed $s,t \in
G$, such that $st^{-1} $ is not in $F$, we see that $(1 \ot
e_F)(\sigma \circ \Delta(\delta_s))(\delta_t \ot 1)=0$, so
$\lambda_s \lambda_{t_1}...\lambda_{t_k} \lambda_t=0$ for any
$t_1,...,t_k \in G.$ Furthermore, the condition that  $\lambda_s
\lambda_{t_1}...\lambda_{t_k} \lambda_t=0$ for any $t_1,...,t_k$
and $t,s$ with $st^{-1}$ not in $F$, actually implies the
apparently stronger condition that  $(\lambda \ot id)(\Delta^n(a))
(\lambda \ot id)(\Delta^n(x_1))...(\lambda \ot id)(\Delta^n(x_k))
  (\lambda \ot id)(\Delta^n(b))=0 \forall x_1,...,x_k \in \cla_0$ whenever
 $a$ and $b$ are such that $(1 \ot e_F)(\sigma \circ \Delta(a))(b
 \ot 1)=0.$ For example, for $a,b$ as above, we have that for any
 $h \in G,$
 $(\lambda \ot id)(\Delta(a))(h).(\lambda \ot
 id)(\Delta(b))(h)=\sum_{s,t} a(sh)b(th) \lambda_s \lambda_t=0$,
 since the above sum is over $(s,t)$ such that
 $sh.(th)^{-1}=st^{-1}$ not in $F$ (for otherwise $a(sh)b(th)=0$),
  and $\lambda_s \lambda_t=0$ for such $s,t.$

  Thus, the $C^*$-algebra $\cle_F$ in this case turns out to be
  the universal  $C^*$-algebra generated by $\{ \lambda_g, g \in G
  \}$ such that each $\lambda_g$ is a nonnegative element,
  $\lambda_s \lambda_{t_1}...\lambda_{t_k} \lambda_{t}=0$ for any
  $t_1,...,t_k $ and $s,t$ with $st^{-1}$ not in $F$, and the
  condition that $\sum_g \lambda_g \lambda_h=\lambda_h$ for all
  fixed $h \in G.$ Then a careful look at \cite{cuntz} tells us
  that the isomorphism conjecture stated by us in the general
  context of discrete quantum groups coincides with the classical
  Baum-Connes conjecture for a discrete group.
\qed

\vspace{4mm}

We now verify the conjecture for any finite dimensional quantum
group.

 \bthm

The conjecture 5.1.4  is true for any finite dimensional quantum
group.

\ethm

{\it Proof :-}\\
Let $\cla$ be a finite dimensional quantum group. We have to prove
that the analytical assembly maps are isomorphism when applied on
the $KK$-groups corresponding to $\cle_F$ for $F=I,$ as $I$ itself
is a finite set in this case. Now, $\cle \equiv \cle_I$ is nothing
but the universal $C^*$-algebra generated by $\{ \lambda(a), a \in
\cla \}$, where $\lambda : \cla \raro \cle $ is CP and unital. To
see this, it is enough to note that the condition (ii) in the
definition of $\cle_I$ is automatic in this case, as $e_I=1$, and
$ (a \ot b) \mapsto \Delta(a)(b \ot 1)$ is a bijection, so that
$\Delta(a)(b \ot 1)=0$ if and only if $a=b=0.$ Furthermore,
$\sum_{\al \in I} \lambda(a)$ is a finite sum in this case, so
this sum actually belongs to $\cle,$ which shows that $\cle$ is
unital.

Now, with the above identification of $\cle$, we claim that $\cle$
is $\cla$-equivariantly homotopic to $\IC$, and hence
$KK^\cla_i(\cle, \IC) \cong KK^\cla_i(\IC,\IC), i=0,1,$ which will
 complete the proof of  the theorem, since we have already proven (\ref{finite}) that
the analytical assembly maps are isomorphisms between
$KK_i^\cla(\IC,\IC)$ and $KK_i(\IC,\hat{\cla}).$ Let us denote the
comultiplication on $\cle$ by $\Delta_\cle.$ To establish our
claim, we first note that since $\cla$ is a finite dimensional
quantum group, there is a unique faithful Haar state on it, say
$\tau,$ and this gives a unital CP map from $\cla$ to $\IC$ given
by $a \mapsto \tau(a).$ From the definition of $\cle,$ there must
exist a unital $C^*$-homomorphsim $\tilde{\tau} : \cle \raro \IC$
such that $\tilde{\tau}(\lambda(a))=\tau(a).$ Let us denote by
$\pi_1$ the homomorphism from $\cle$ to itself given by
$\pi_1(a)=\tilde{\tau}(a)1_\cle,$ and by $\pi_0$ the identity map
from $\cle$ to itself. To verify our claim, we need to give an
equivariant homotopy connecting $\pi_0$ and $\pi_1$. Define for $t
\in [0,1]$ the unital CP map $\lambda_t : \cla \raro \cle$ given
by $\lambda_t(a): (1-t)\lambda(a)+t \tau(a)1_\cle.$ By the
universality of $\cle$, there exist unital $C^*$-homomorphisms
$\pi_t : \cle \raro \cle$, satisfying $\pi_t \circ \lambda =
\lambda_t.$ Clearly, $t \mapsto \pi_t(x)$ is continuous for $x$
belonging to the dense $\ast$-algebra spanned by the elements in
the range of $\lambda$, and hence for all $x.$ It remains to check
that $\pi_t$ are equivariant. To this end, it is enough to show
that $(\pi_t \ot
id)(\Delta_\cle(\lambda(a)))=\Delta_\cle(\pi_t(\lambda(a)))$ for
all $a \in \cla.$ Indeed, we have,

\bean
\lefteqn{(\pi_t \ot id)(\Delta_\cle(\lambda(a))}\\
 &=& \Delta_\cle(\lambda_t(a))\\
 &=& (1-t) \Delta_\cle(\lambda(a))+t \tau(a)1\\
 &=& (1-t)(\lambda \ot id)(\Delta(a))+t(\tau \ot id)(\Delta(a))\\
 &=& (\lambda_t \ot id)(\Delta(a))\\
 &=& (\pi_t \ot id)((\lambda \ot id)(\Delta(a)))\\
 &=& (\pi_t \ot id)(\Delta_\cle(\lambda(a))),\\
\eean
 which completes the proof. Note that we have used in the above
 the fact that $\tau(a)1=(\tau \ot id)(\Delta(a)),$ which follows
 from the definition of the Haar state.
 \qed

\brmrk

We have seen in Example 2 of section 4 that $\mu_0$ is surjective
from $KK_0^\cla(\cla,\IC)$ to $KK_0(\IC,\hat{\cla})$ for
$\cla=\hat{SU_q(2)}.$ From this, it is possible to argue that
$\mu_0^r :lim_F KK_0^\cla(\clc,\IC) \raro KK_0(\IC,\hat{\cla}_r)$
is surjective. Indeed, the limit algebra $\cle \equiv lim_F
\cle_F$ is nothing but the universal $C^*$ algebra generated by
$\lambda(a), a \in \cla$ with $\lambda$ CP and unital from
 $\clm(\cla)$ to $\clm(\cle).$ From the universality, one can get
 an equivariant  $C^*$-homomprhism from $\cle$ to $\cla$, which in
 turn gives a homomorphism from the equivariant KK-group of $\cla$
 to  $lim_F KK_0^\cla (\cle_F, \IC)$ and the composition of the analytical assembly
 map with this will coincide with the analytical assembly map from
 $KK_0^\cla(\cla,\IC)$ to $KK_0(\IC,\hat{\cla}).$ From this,
 surjectivity will follow.

 \ermrk

 We postpone the
attempts to verify this conjecture for some non-classical quantum
groups. e.g. duals of $SU_q(n)$'s, for future works.

{\it Acknowledgement :}\\
 D. Goswami  would like to express his gratitude to I.C.T.P. for a visting
research fellowship during January-August 2002, and to  A.O. Kuku
and  the other organisers of the ``School and Conference on
Algebraic K Theory and Its Applications" at I.C.T.P. (Trieste) in
July 2002. He would also like to thank
  T. Schick  for sending some relevant preprint,  and   A.
Valette and  I. Chatterjee for giving useful information regarding
 some  manuscript (yet to be published)  by A. Valette.

\end{document}